\documentclass[11pt,thmsa,reqno]{amsart}

\usepackage[all]{xy}
\usepackage{color}
\usepackage{mathtools} 
\usepackage{fancybox}
\usepackage{amssymb}
\usepackage[hidelinks]{hyperref} 
\usepackage{comment}
\usepackage{graphicx}
\usepackage{faktor}
\usepackage[normalem]{ulem}
\usepackage{tikz-cd}

 \newtheorem{theorem}{Theorem}[section]
 \newtheorem*{thm*}{Theorem} 
 \newtheorem{lemma}[theorem]{Lemma}
 \newtheorem{prop}[theorem]{Proposition}
 
 \newtheorem{corollary}[theorem]{Corollary}
 
  \theoremstyle{definition}
 \newtheorem{example}[theorem]{Example}
 \newtheorem{remark}[theorem]{Remark}

\renewcommand{\graph}{\mathrm{graph}}

\newcommand{\id}{\mathrm{id}}

\newcommand{\vX}{\mathfrak{X}}

\newcommand{\blupC}{\mathrm{Blup}}
\newcommand{\blupiotaA}{\mathrm{Blup}}

\newcommand{\cF}{\mathcal{F}}

\newcommand{\cL}{\mathcal{L}}

\newcommand{\RR}{\ensuremath{\mathbb R}}

\newcommand{\pd}[1]{\frac{\partial}{\partial #1}} 

\renewcommand{\graph}{\mathrm{graph}}
\newcommand{\B}{\mathfrak{B}}

\newcommand{\pr}{\mathrm{pr}}
\newcommand{\rank}{\mathrm{rank}}

\newcommand{\blup}{\mathrm{Blup}}
\newcommand{\D}{\mathrm{d}}

\begin{document}

 \author{Andreas Sch\"{u}\ss ler and Marco Zambon}
\address{KU Leuven, Department of Mathematics, Celestijnenlaan 200B box 2400, 3001 Leuven, Belgium.}

\email{marco.zambon@kuleuven.be}
\address{KU Leuven, Department of Mathematics, Celestijnenlaan 200B box 2400, 3001 Leuven, Belgium.}
   
\subjclass[2020]{Primary: 53D17. Secondary: 57R30 \\Keywords: Lie algebroid, singular foliation, Dirac structure, pullback, blow-up.}

 \title[A note on pullbacks and blowups]{A note on pullbacks and blowups of Lie algebroids, singular foliations, and Dirac structures}

\begin{abstract}
Lie algebroids, singular foliations, and Dirac structures are closely related objects. We examine the relation between their pullbacks under maps satisfying a constant rank or transversality assumption. A special case is given by blowdown maps. In that case, we also establish the relation between the blowup of a Lie algebroid and its singular foliation. 
\end{abstract}

\maketitle



   \section{Introduction}
\label{sec:intro}
{
Dirac structures are a kind of geometric structure on manifolds that generalizes both closed 2-forms and Poisson structures. {In the first part of this note we} determine the singular foliation underlying the pullback of a Dirac structure $L$ under a {suitable} map, by showing that it coincides with the pullback of the singular foliation of $L$. Here ``singular foliation'' is understood as a module of vector fields, in the sense of \cite{AndrSk}, and not just as the underlying partition of the manifold into leaves. We do this in Corollary \ref{cor:Diracsf}. 
}

{To obtain Corollary \ref{cor:Diracsf}, we involve Lie algebroids.} 
It is well-known that
any Lie algebroid $A$ induces a singular foliation ${\cF_A}:=\rho(\Gamma_c(A))$, where $\rho$ denotes the anchor. Further, any Dirac structure inherits a Lie algebroid structure, with anchor the restriction of $\pr_{TM}$, and bracket the restriction of the Courant bracket.
{Each of these three structures--Dirac structures $L$, Lie algebroids $A$, singular foliations $\cF$--can be pulled back along smooth a map $f$ satisfying compatibility conditions. We denote the pulled-back structures by 
$\mathfrak{B} L$, $f^!A$, and $f^{-1}\cF$, respectively.
The pullback of the Lie algebroid underlying a Dirac structure $L$ will be denoted by $f^!L$.
}

Given a smooth map $f\colon B\to M$ and a Lie algebroid $A$ over $M$, we say that ``$f^!A$ {is smooth}'' if the subspaces $\rho(A_{f(x)})+(f_*)(T_xB)$ have the same rank for
 all $x\in B$. This condition is equivalent to the  existence of the pullback Lie algebroid {$f^!A$} of $A$ by $f$. {If
$A$ is the Lie algebroid underlying a Dirac structure $L$, by 
 \cite[Thm. 7.33]{PoisGeoBookAMS} this condition also ensures that the pullback Dirac structure exists, i.e.\ that $\mathfrak{B} L$ is a Dirac structure.}

{To state the main results of the first part of this note, we paraphrase Proposition \ref{prop:LAfol},  Proposition \ref{prop:nottrans} together with  Proposition \ref{prop:iso}, and Corollary \ref{cor:Diracsf}.
\begin{itemize}
    \item Let $A$ be a Lie algebroid on $M$, assume that $f^!A$ is smooth. Then   $$\cF_{f^!A}=f^{-1}(\cF_A).$$
In particular, $f^{-1}(\cF_A)$
    is a singular foliation.
    \item Let $L$ be a Dirac structure on $M$, assume that $f^!L$ is smooth,  Then 
$\mathfrak{B} L$ is a Dirac structure, and $$\cF_{\mathfrak{B} L}=\cF_{f^!L}.$$ 
  If we assume the stronger condition that $f$ is transverse to $L$, then there is a canonical Lie algebroid isomorphism $$\mathfrak{B} L\cong f^!L.$$
\item Let $L$ be a Dirac structure on $M$, assume that $f^!L$ is smooth. Then  
$\mathfrak{B} L$ is a Dirac structure, and  $$\cF_{\mathfrak{B} L}= f^{-1}(\cF_L).$$
\end{itemize}
}

{The third item, Corollary \ref{cor:Diracsf}, is obtained combining the first and second item.}
We emphasize that Corollary \ref{cor:Diracsf} is in a sense  an optimal result, as we explain just before Example \ref{ex:counter}.

\bigskip
In the second part of this note 
we consider a submanifold $N\subset M$ and the blowup $B$ of $M$ along $N$. The blowdown map $p\colon B\to M$ is a diffeomorphism on an open dense subset. If a Dirac structure $L$ on $M$ is transverse to $N$, it admits a unique lift to $B$. Using the above results we show that the singular foliation of the lift is the pullback of the singular foliation of $L$.

If a Lie algebroid $A$ on $M$ is transverse to $N$, then it gives rise to several Lie algebroids on $B$: one of them is the pullback Lie algebroid $p^!A$, others are  obtained blowing up $A$ itself w.r.t.\ a Lie subalgebroid  $C$ supported on $N$. 
 We determine the singular foliation {$\cF_{\blupC}$} of the blown-up Lie algebroid for any Lie subalgebroid $C$ over $N$ which contains the isotropies of $A$ over $N$, and make the result more explicit in two cases:
 \begin{itemize}
    \item when $C$ is the restriction of $A$ to $N$, {$\cF_{\blupC}$} is given by the intersection of the pullback of the singular foliation of $A$ with the $b$-tangent bundle of $B$ with respect to the hypersurface $p^{-1}(N)$ (see
 Example \ref{ex:restricted_LA}), 
 \item when $C$ is the isotropy Lie algebroid of $A$ over $N$, {$\cF_{\blupC}$} is given by the intersection of the pullback of the singular foliation of $A$ with the edge Lie algebroid \cite{fineedge} associated to the fibration $p^{-1}(N)\to N$  (see
 Example \ref{ex:isotropy_LA}).
\end{itemize}

\bigskip
{
\noindent\textbf{Relation with the literature:} 
A special case of Proposition \ref{prop:LAfol}, in which the map $f$ is assumed to be a surjective submersion, appeared in \cite[Lemma 1.15]{HME}.

The statement of Proposition \ref{prop:iso} appears in
\cite[\S 5.1]{SplittingThmEulerlike}, but no further details are given there. A proof can be obtained from the one of Proposition 6.6 in
\cite[\S 6.2]{meinrenkenPoisNotes}, yielding an isomorphism which is the inverse of the one we construct in  our proof  of Proposition \ref{prop:iso}.

\bigskip
\noindent {\noindent\textbf{Notation:}  Given a  Lie algebroid $A$, we   denote the corresponding singular foliation by $\cF_A:=\rho(\Gamma_c(A))$, where $\rho$ is the anchor. Further, given a  Dirac structure $L$,
we denote  by $\cF_L$ the singular foliation of the  underlying Lie algebroid 
(hence $\cF_L=\pr_{TM}(\Gamma_c(L))$).
}

\bigskip
\noindent\textbf{Acknowledgements:} 
{M.Z. thanks Henrique Bursztyn for his feedback on the relation with the literature.}
{We are grateful to the referee for her/his insights, which allowed to simplify considerably some arguments, improve some of the results, and obtain new ones.}
We acknowledge partial support by Methusalem grant METH/21/03 - long term structural funding of the Flemish Government.
M.Z. acknowledges partial support by EOS project G0I2222N and by FWO project G0B3523N  (Belgium).

 \section{Definitions of pullbacks}
\label{sec:def}

Let $f\colon B\to M$ be a smooth map.

\subsection{Pullbacks of Lie algebroids}
\label{subsec:PullLA}
A \emph{Lie algebroid} is a vector bundle $A$ over a manifold $M$, together with a bundle map $\rho\colon A\to TM$ called anchor and a Lie bracket on the smooth sections $\Gamma(A)$, such that $[a,fa']=\rho(a)(f) a'+f[a,a']$ for all sections $a,a'$ and functions $f$. Lie algebroids play an important role in differential geometry, generalizing both tangent bundles and Lie algebras.

Let $A$ be a Lie algebroid over $M$. 
We denote by $f^*A\to B$ the pullback of $A$ as a vector bundle. 
{Consider the following subset of the vector bundle $f^*A\oplus TB$:
$$f^!A:=f^*A\times_{TM}TB,$$ the fiber product over $\rho\colon A\to TM$ and $f_*\colon TB\to TM$. 
Assume that \emph{$f^!A$ {is smooth}}, meaning that 
the fibers of the map $f^!A\to B$ have constant rank.
{This assumption is equivalent to}  $f^!A$ being a vector subbundle of $f^*A\oplus TB$, as we recall in Lemma \ref{rem:f^!Aexists}.}

We call $f^!A$ 
the \emph{pullback of $A$ as a Lie algebroid} (\cite{higgins.mackenzie:1990a}, see also \cite[\S 4]{MacKbook}). 
Any section of $f^!A$ is of the form\footnote{Here we slightly abuse notation, writing $a_i$ instead of $a_i\circ f$.} $(\sum h_i a_i, X)$ for finitely many $h_i\in C^{\infty}(B)$ and $a_i\in \Gamma(A)$, and for $X\in \mathfrak{X}(B)$ such that $f_*X=\sum h_i \rho(a_i)$. Here, we view $f_*$ as a map $f_*\colon TB\to f^*TM$ covering $\id_B$.

More importantly, $f^!A$ is a Lie algebroid, with the second projection as anchor. The bracket of two sections as above is determined by the Lie algebroid bracket of $A$ and the Leibniz rule:
$$
 \Big[  \Big( \sum_i h_i a_i, X\Big),\Big( \sum_j h'_j a_j', X'\Big)\Big]=(\star, [X,X']),
 $$
 where
 \begin{equation}\label{eq:bracket}
 \star= \sum_{i,j} h_i h'_j [a_i,a_j']+\sum_j X(h'_j)a_j'-\sum_i X'(h_i)a_i.
\end{equation}
{Upon the identification $B\cong \graph(f)\subset M\times B$, there is an identification of  $f^!A$ with  
the preimage of $T\graph(f)$ under the anchor  in
the product Lie algebroid $A\times TB$, which indeed is a Lie subalgebroid over $\graph(f)$, see \cite[\S 7.4]{meinrenkenGroidsAlgoidsNotes}, \cite{meinrenken2024liealgebroids}.}

 Further, the first projection yields a Lie algebroid morphism  $f^!A\to A$  \cite[\S 4.3]{MacKbook}. This morphism is not surjective in general; for instance, for any smooth map $f\colon B\to M$, we have $f^!TM=TB$, and the morphism is given by the derivative $f_*\colon TB\to TM$. 

\begin{lemma}\label{rem:f^!Aexists}
{Let $f\colon B\to M$ be a smooth map, and $A$ a Lie algebroid over $M$. Then $f^!A$  {is smooth} if{f} the following subspaces have constant rank for all $x\in B$:
\begin{equation}\label{eq:transverse}
\rho(A_{f(x)})+(f_*)(T_xB).
\end{equation}
{In turn, this is equivalent to  $f^!A$ being a} vector subbundle of $f^*A\oplus TB$.}  
\end{lemma}

\begin{proof}
 {Recall that the meaning of ``$f^!A$ is smooth'' is simply that 
the fibers of the map $f^!A\to B$ have constant rank.
 Notice that
$f^!A=\ker(\phi)$, where $\phi$ is the vector bundle map 
$$
\phi\colon f^*A\oplus TB\to f^*TM,\; (a,v)\mapsto \rho(a)-f_*v.   
$$
By the rank-nullity theorem, the fibers of $\ker(\phi)$ have constant rank if{f} the fibers of the image $\text{im}(\phi)$ have constant rank. The latter are given precisely by 
\eqref{eq:transverse}.
}

{
In this case, $f^!A$ is the kernel of a constant rank map defined on $f^*A\oplus TB$ (namely $\phi$), thus it is a vector subbundle. Conversely, if $f^!A$ is a  a vector subbundle, its fibers obviously have constant rank. This proves  the second assertion. 
}
\end{proof}

{A special case of Lemma \ref{rem:f^!Aexists} occurs when $A$ is transverse to $f$, in the sense that the subspaces \eqref{eq:transverse} equal the whole of 
$T_{f(x)}M$.  
 In that case, 
$f^!A$
has rank equal to $\rank(A)+(\dim B -\dim M)$.
Indeed, this condition implies that the map $\phi$ is transverse to the zero section.}

\subsection{Pullbacks of singular foliations}
A \emph{singular foliation} $\cF\subset \vX_c(M)$ is a $C^{\infty}(M)$-submodule of the compactly supported vector fields, which is locally finitely generated and involutive w.r.t.\ the Lie bracket \cite{AndrSk}.
Here, we use the subscript $c$ to denote ``compactly supported''.
A singular foliation gives rise to a decomposition of $M$ into immersed submanifolds of possibly varying dimension, called leaves. {In general, the singular foliation can not be recovered from the decomposition into leaves. For instance, on $M=\RR$, the singular foliations $C^{\infty}_c(M) x^k\pd{x}$ are distinct for all integers $k\ge 1$, but they all have the same underlying partition into leaves (namely the origin, the positive, and the negative axis).
}

Suppose we have a singular foliation $\cF\subset \vX_c(M)$ which is transverse to $f$, i.e.\ each leaf of the singular foliation is transverse to $f$. 
Then there exists a \emph{pullback singular foliation} \cite[\S 1.2.3]{AndrSk} given by
\begin{equation*}
f^{-1}\cF:=\{X\in \vX_c(B): f_*X=\sum h_i (Y_i\circ f) \text{ for } h_i\in C^{\infty}_c(B), Y_i\in \cF\},
\end{equation*}
where the sum is finite. The leaves of $f^{-1}\cF$ are the connected components of the preimages of the leaves of $\cF$. {If $\cF$ is not transverse to $f$, then 
$f^{-1}\cF$ is an involutive submodule which might fail to be locally finitely generated, see Example \ref{ex:fail} below.}

\subsection{Pullbacks of Dirac structures}
A \emph{Dirac structure}  $L\subset TM\oplus T^*M$ on $M$ is a maximal isotropic subbundle (w.r.t.\ the canonical symmetric pairing) that is involutive w.r.t.\ the Courant bracket 
\begin{equation}\label{eq:courantbracket}
[(X,\xi),  (X',\xi')]=([X,X'], \cL_X \xi'-\iota_{X'} \D \xi).
\end{equation}
Prototypical examples of Dirac structures are graphs of closed 2-forms and of Poisson bivector fields.

Consider 
$$\mathfrak{B}L:=\{(X,f^*\eta):(f_*X,\eta)\in L\},$$
a collection of maximal isotropic subspaces of $TB\oplus T^*B$.
Assume that  $\rho(L)$ is transverse to $f$, i.e.\ {the subspaces \eqref{eq:transverse} equal $T_{f(x)}M$} for $L=A$, $\rho=\pr_{TM}$. Then  $\mathfrak{B}L$ is a Dirac structure on $B$ \cite[Proposition 5.6.]{Bursztyn_2013}.
The map $f$ from $(B, \mathfrak{B}L)$ to $(M,L)$ is then said to be a \emph{backward Dirac map}. {Without the transversality assumption, the (constant rank)  collection of maximal isotropic subspaces  
$\mathfrak{B}L$ can fail to be a smooth subbundle, see Example \ref{ex:fail} below.}

\subsection{Examples}

{Let $f$ be a smooth map. 
We saw above that, given  a singular foliation or Dirac structure to which  $f$ is transverse, or  a Lie algebroid $A$  such that $f^!A$  {is smooth},  we can always pull back this structure. 
The transversality condition is guaranteed when $f$ is a submersion (a condition on $f$ alone).
If $A$ is an involutive distribution  on $M$ (i.e. a distribution tangent to a  regular foliation), and $\iota\colon N\hookrightarrow M$ a submanifold  such that $A_{x}+T_xN$ has constant rank for all $x\in N$, then $\iota^!A$ {is smooth}, and is indeed an involutive distribution on $N$.}

{
We now display an example to show that, without suitable  assumptions, 
pullbacks might not exist.}

\begin{example}
\label{ex:fail}
{
Consider a smooth\footnote{Notice that  $\varphi$ can not be an analytic function.} function $\varphi \colon \RR\to \RR$ such that $\varphi|_{\RR_{\ge 0}}=0$ and 
$\varphi|_{\RR_{< 0}}$ is nowhere vanishing, with nowhere vanishing derivative. Consider the submanifold $N=\graph(\varphi)$ of  $M=\RR^2$. On $M$
consider the Dirac structure $$L=\RR \pd{x}\oplus \RR dy,$$ the underlying Lie algebroid $A$, and the underlying singular foliation $\cF$ generated by $\pd{x}$. Notice that $N$ is transverse to these structures only at points $(x,\varphi(x))$ with $x<0$.
We have: 
\begin{itemize}
    \item $\iota^!A$ {is not smooth}. 
Indeed,   $\iota^!A$ has rank $1$ at points $(x,\varphi(x))$ with $x<0$, and rank $2$ where $x\ge 0$. 
\item The involutive submodule $$\iota^{-1}\cF= \left\{(\Phi^* h)\pd{x}: h\in C^{\infty}_c(\RR) \text{ satisfies } h |_{\RR_{< 0}}=0\right \}$$ 
is not locally finitely generated. Here, $\Phi\colon N\to \RR$ is the natural diffeomorphism given by the first projection.
\item  $\mathfrak{B} L$ is not a smooth subbundle of $TN\oplus T^*N$. Indeed, it equals $\{0\}\oplus T^*N$ at points $(x,\varphi(x))$ with $x<0$, and it equals  $TN\oplus \{0\}$ for $x\ge 0$.
\end{itemize}
}
\end{example}

 \section{Lie algebroids and singular foliations}
\label{sec:la}

Under {suitable} assumptions, 
the operations of taking the singular foliation of a Lie algebroid and of taking the pullback commute. 

\begin{prop}\label{prop:LAfol}
  Let $f\colon B\to M$ be a smooth map.
Let $A$ be a Lie algebroid over $M$, 
and assume that $f^!A$ {is smooth}. Then
$${\cF_{f^!A}=f^{-1}(\cF_A).}$$
\end{prop}

\begin{proof}
``$\subset$'' Any compactly supported section of $f^!A$ is of the form $(\sum h_i a_i, X)$ for finitely many $h_i\in C^{\infty}_c(B)$ and $a_i\in \Gamma(A)$, and $X\in \mathfrak{X}_c(B)$ such that $f_*X=\sum h_i \rho(a_i)$. We may assume that the $a_i$ are compactly supported, by multiplying them with a compactly supported function which equals $1$ on $f(\mathrm{Supp
}(h_i))$.
Since $\rho(a_i)\in {\cF_A}$, the conclusion follows.

``$\supset$''  Let $X\in f^{-1}({\cF_A})$, i.e.\ $f_*X=\sum h_i (Y_i\circ f)$, where $h_i\in C^{\infty}_c(B), Y_i\in {\cF_A}$. We have $Y_i=\rho(a_i)$ for certain $a_i\in \Gamma_c(A)$. Then $(\sum h_i a_i, X)\in \Gamma_c(f^!A)$, since $\sum h_i \rho(a_i)=\sum h_i Y_i=f_*X$. 
This is an element of $\Gamma_c(f^!A)$ whose image under the anchor is $X$.
\end{proof}

{Notice that Proposition \ref{prop:LAfol} implies in particular that $f^{-1}(\cF_A)$  is a singular foliation, providing a criterion (different from transversality  \cite[Prop. 1.10]{AndrSk}) to ensure that the pullback of a singular foliation is again  a singular foliation.}

 \section{Dirac structures and Lie algebroids}
\label{sec:Dirac}

  Let $L$ be a Dirac structure over $M$, 
  {and $f\colon B\to M$   a smooth map.}

\begin{prop}\label{prop:iso}
  Let $L$ be a Dirac structure on $M$, transverse to $f$.
Consider  $\mathfrak{B} L$, viewed as a Lie algebroid, and the Lie algebroid pullback $f^!L$ of $L$. 
{There is a canonical isomorphism (covering $\id_B$) of Lie algebroids  $$\mathfrak{B} L\cong f^!L.$$}
\end{prop}
\begin{proof} 
Notice first that at a point $x\in B$, any element of $(f^!L)_x$ is of the form $(\ell,X)$ where $\ell\in L_{f(x)}$, $X\in T_xB$ with $f_*X=\pr_{TM}\ell$, 
hence $\ell=(f_*X,\eta)$ for some $\eta\in T^*_{f(x)}M$.
Consider the vector bundle map
\begin{equation}\label{eq:psi}
   \psi\colon f^!L \to \B L,\;\; ((f_*X,\eta),X)\mapsto (X,f^*\eta). 
\end{equation}

This vector bundle map clearly takes values in   $\B L$ .
We first argue that $\psi$ is a vector bundle isomorphism.
It is surjective, by the very definition of $\B L$. We have $\rank (f^!L)=\dim B$ using $\rank (L)=\dim M$ {by the transversality assumption, see the text after Lemma \ref{rem:f^!Aexists}}.
Hence, $\psi$ is a map between vector bundles of the same rank. Therefore, it is also injective, and an isomorphism.

The map $\psi$ clearly preserves the anchors. To show that it preserves brackets, {one can take sections $\sigma$ of $\B L\subset TB\oplus T^*B$ and $\tau$ of $L\subset TM\oplus T^*M$ which are related by $f$, meaning that 
$\sigma=(Z,f^*\eta)$ and 
$\tau=(f_*Z,\eta)$ for an $f$-projectable vector field $Z\in \vX(B)$ and for $\eta\in \Omega^1(M)$; the statement then follows from the fact \cite[Lemma 6.1]{meinrenkenPoisNotes} that the Courant bracket of related sections are again related.
}

{Alternatively, one can show that $\psi$ preserves brackets by a direct computation:} 
take a section $$(\sum h_i (Y_i, \eta_i), X)\in \Gamma(f^!L),$$  where $h_i\in C^{\infty}(B)$ and $(Y_i, \eta_i)\in \Gamma(L)$, for $X\in \mathfrak{X}(B)$ such that $f_*X=\sum h_i Y_i$. Under $\psi$, it is mapped to $(X, \sum h_i  f^* \eta_i)\in \Gamma(\B L)$. Use eq.\ \eqref{eq:bracket} to compute the bracket of two sections of $\Gamma(f^!L)$, {recalling that the Lie bracket of $L$ is the restriction of the Courant bracket on $M$.} Use the Courant bracket \eqref{eq:courantbracket} on $B$ to compute the bracket of their images under $\psi$, together with identities such as $\iota_X f^*\eta=\sum_i h_i \iota_{Y_i} \eta$ and  the fact that $L$ is isotropic, to conclude that $\psi$ preserves brackets.
\end{proof}

\begin{remark}
  We provide a direct proof that the map $\psi$ in \eqref{eq:psi} is injective, without using any dimension considerations.
  
We first claim: \emph{The vector bundle map over $\id_B$ $$\phi \colon f^*L\to f^*TM\oplus T^*B, \;\; (Y,\eta)\mapsto (Y, f^*\eta)$$ is injective.}

Indeed, for all $x\in B$, the tranversality condition \eqref{eq:transverse} (for $A=L)$ can be rephrased as $L_{f(x)}\cap [f_*(T_xB)]^{\circ}=\{0\}$, as one sees taking annihilators and using $\rho(L)=(L\cap TM)^{\circ}$. Now let $(Y,\eta)\in (f^*L)_x$ lie in the kernel of $\phi$, i.e.\ $Y=0$ and $\eta \in [f_*(T_xB)]^{\circ}$. Then $(Y,\eta)\in L_{f(x)}\cap [f_*(T_xB)]^{\circ}$, so it vanishes. This proves the claim.

Let $((f_*X,\eta),X)\in f^!L$ lying in the kernel of $\psi$. Then $(X,f^*\eta)=0$, and in particular $(f_*X,f^*\eta)=0$.
But $(f_*X,f^*\eta)$ is the image of $(f_*X,\eta)\in f^*L$ under the map $\phi$ above. The injectivity of $\phi$ implies that $(f_*X,\eta)=0$. Hence, $\psi$ is injective.
\end{remark}

We now consider a more general setting than the one of
Proposition \ref{prop:iso}, by replacing the transversality assumption  with the  weaker requirement that 
$f^!L$  is smooth. {In the following proposition, the result that $\mathfrak{B} L$ is a Dirac structure is due to \cite[Thm. 7.33]{PoisGeoBookAMS}; we state again the result and its proof, since we need them for the second part of the proposition.}

\begin{prop}\label{prop:nottrans} 
{Let $f\colon B\to M$ be  a smooth map, and   $L$   a Dirac structure on $M$ such that $f^!L$  is smooth.}  Then 
$\mathfrak{B} L$ is a Dirac structure. Further, $\mathfrak{B} L$ and $f^!L$   induce the same singular foliation.
\end{prop}
\begin{proof}
{
Since   $f^!L$ is a vector  subbundle of $f^*L \oplus TB$,
(see Lemma \ref{rem:f^!Aexists}),
we can view the map $\psi$ in \eqref{eq:psi} as a (smooth) vector bundle map $\Psi\colon f^!L\to TB\oplus T^*B$. The image of $\Psi$ is $\mathfrak{B} L$, which pointwise  consists of maximal isotropic subspaces. Hence $\Psi$ has constant rank, and therefore, as the base map of $\Psi$ is the identity, its image $\mathfrak{B} L$ is a smooth subbundle of $TB\oplus T^*B$. The involutivity follows from the proof of  Proposition \ref{prop:iso}.}

{Since the vector bundle map $\psi$ in \eqref{eq:psi} commutes with the projections to $TB$, the singular foliation induced by $f^!L$ is contained in the one induced by $\mathfrak{B} L$. Since we have established that $\psi$   is surjective, every section of $\mathfrak{B} L$ is the image of a section of $f^!L$, implying that
the two singular foliations agree.
}
\end{proof}

\begin{remark}
In the set-up of Proposition \ref{prop:nottrans}, $f^!L$ can have strictly larger rank than $\mathfrak{B} L$. Consider for instance the inclusion $\iota$ of a point $p$ in $M$. For any Dirac structure $L$ we have $\mathfrak{B} L=\{0\}$, while  
the pullback of any  Lie algebroid  is its isotropy Lie algebra $\ker \rho_p$. 
    \end{remark}

The following corollary\footnote{This corollary 
is stronger than the analog statement for leaves, which is certainly known {in the transverse case}, and which states the following: the leaves of {$\cF_{\mathfrak{B} L}$} are the connected components of the preimages of the leaves of {$\cF_L$}.} follows immediately from Propositions \ref{prop:LAfol} and {Proposition \ref{prop:nottrans}}.

\begin{corollary}\label{cor:Diracsf}
   Let $L$ be a Dirac structure over $M$ such that $f^!L$ {is smooth}.
   {Then  
$\mathfrak{B} L$ is a Dirac structure, and  $$\cF_{\mathfrak{B} L}= f^{-1}(\cF_L).$$}
\end{corollary}

{There do exist cases in which $\mathfrak{B} L$ is a 
  (smooth) Dirac structure but $f^!L$ {is not smooth}. The conclusion of  Corollary \ref{cor:Diracsf} does \emph{not} hold in general if we only assume that $\mathfrak{B} L$ is a Dirac structure. We give a counterexample in item ii) below, following \cite{Coregularsubmanifolds}}. {A similar counterexample is obtained also following \cite[Rem. 7.35]{PoisGeoBookAMS}.}

  \begin{example}\label{ex:counter}
{  i) Let $M=\mathfrak{so}(3)^*$, with  Dirac structure $L$ given by the graph of the canonical linear Poisson structure. Let $f\colon B\to M$ the blowdown map defined on the blowup $B$ of $M$ at the origin. Then $\mathfrak{B} L$ is a 
  (smooth) Dirac structure   \cite[Thm. 4.2]{PMCT3}, \cite[Ex. 7.1]{MSZ}. However, $f^!L$ {is not smooth}:  the l.h.s.\ of \eqref{eq:transverse} has rank $3$ away from the exceptional divisor $f^{-1}(0)$ (since $f$ defines a diffeomorphism there), but has rank $1$ at points of the divisor. Nevertheless, one can check that   {$\cF_{\mathfrak{B} L}= f^{-1}(\cF_L)$}.
  }  

{ii) We revisit \cite[Ex. 2.7]{Coregularsubmanifolds}, which the authors use to exhibit a pathology they call ``jumping phenomenon'' (and remark that this pathology does not arise for the class of coregular submanifolds of Poisson manifolds).
Fix $f\in C^{\infty}(\RR^2)$, and 
consider the embedding $$\iota\colon  \RR^2\to \RR^4,\;\;(x,y)\mapsto (x,y, f(x,y)^2, f(x,y)^2).$$
Endow $\RR^4$ with the Dirac structure $L$ given by the graph of the Poisson structure $\pi=\pd{x_1}\wedge \pd{x_2}+x_3\pd{x_3}\wedge \pd{x_4}$.
The authors check that $\mathfrak{B} L$ is the (smooth) Dirac structure $\graph(\pd{x}\wedge \pd{y})$, whose underlying singular foliation is 
${\cF_{\mathfrak{B} L}}=\vX_c(\RR^2)$. They point out that its (unique) leaf may not be contained in any leaf of the ambient Poisson manifold.
}

{We check that $\iota^{-1}{(\cF_L)}$ is strictly contained in ${\cF_{\mathfrak{B} L}}$, thus showing that the conclusion of Corollary \ref{cor:Diracsf} does not hold in this case. Here,
{$\cF_L$} is the singular foliation induced by the Dirac structure $L$, which is  generated by $\pd{x_1}, \pd{x_2},x_3\pd{x_3},x_3\pd{x_4}$. To do so, we check that the product of $\pd{x}$ with any compactly supported function is {generally} not an element of $\iota^{-1}{(\cF_L)}$.
Indeed, denoting $N:=\iota(\RR^2)$, we have
\begin{equation}\label{eq:ix}
    \iota_*\pd{x}=\left(\pd{x_1}+2f f_{x_1} \pd{x_3}+2f f_{x_1}\pd{x_4}\right)|_N,\end{equation}
where $f$ is viewed as a function of $(x_1,x_2)$ and $f_{x_1}$ denotes its first partial derivative.
In general, \eqref{eq:ix} can not be written as a linear combination
$$\left(\pd{x_1}+h_3 x_3 \pd{x_3}+h_4 x_3 \pd{x_4}\right)|_N$$
for $h_3,h_4\in C^{\infty}(N)$.
To see this,  compare the coefficients of $\pd{x_3}$:  on $N$ we have $x_3=f(x_1,x_2)^2$, and in general the equation $$2 f_{x_1}=h_3 f$$ does not have a smooth solution $h_3$ (take for instance $f(x_1,x_2)=x_1$).
}
   \end{example}

 \section{Blowups}
\label{sec:blowup}

Let $N$ be a closed and embedded submanifold of a manifold $M$, denote by $$B:=\blup(M,N)$$ the  \emph{real projective blowup} along $N$ \cite{griffith.harris:1978a}. As a set, $B$ is the disjoint union of $M\setminus N$ and the projectivization $\mathbb{P}$ of the normal bundle $TM|_N/TN$ of $N$. It comes with a smooth and proper map $p\colon B\to M$, called \emph{blowdown map}, which restricts to a diffeomorphism from $B\setminus \mathbb{P}$ to  $M\setminus N$.
The codimension $1$ submanifold $\mathbb{P}=p^{-1}(N)$ is called \emph{exceptional divisor}. 
Every vector field tangent to $N$ admits a (unique) lift to a vector field $\widetilde{Y}$ on $B$ which is $p$-related to $Y$, see e.g.\ \cite[Proposition 1.5.40]{singfolnotes}, \cite[Lemma 3.5]{schuessler:2024a}.

\begin{example}
    The blowup is particularly easy to describe when $M$ is the total space of a vector bundle $E\to N$. In that case, $\blup(E,N)$ is the tautological line bundle over the projectivization of $E$. For instance, $\blup(\mathbb{\mathbb{\RR}}^2,\{0\})$ is the non-trivial line bundle over $\mathbb{RP}^1=S^1$, i.e.\ the M\"obius strip.
\end{example}

\subsection{Lifted Dirac structures on the blowup}
Let $L$ be a Dirac structure over $M$ for which $N$ is a transversal, i.e.\    
 $\rho(L_y)+T_yN=T_yM$ for all $y\in N$. 
At every   $x\in p^{-1}(N)$, the derivative of the blowdown map $p$ satisfies
$T_{p(x)}N\subset p_*(T_xB)$. Hence, the map $p$ is transverse to $L$, and the assumptions of Corollary \ref{cor:Diracsf} are satisfied, yielding:

\begin{corollary}
{Let $L$ be a Dirac structure over $M$ for which $N$ is a transversal. Let $p$ be the blowdown map. 
Then $\mathfrak{B} L$ is a Dirac structure and $$\cF_{\mathfrak{B} L}=p^{-1}(\cF_L).$$}
\end{corollary}

\subsection{Blowups and Lie algebroids}

Let $\pi\colon A\to M$ be a Lie algebroid over $M$. Via the blowdown map, one can lift $A$ to a Lie algebroid structure on $B\setminus \mathbb{P}$. In general, there are several distinct\footnote{This is in contrast to the case of Dirac structures, where the uniqueness is forced by the fact that Dirac structures over $B$ are subbundles  
of a prescribed vector bundle, namely $TB\oplus T^*B$.} extensions to Lie algebroids over the whole of $B$. In this subsection, we want to describe the singular foliations of  some  of them.

Some possible extensions   
are given by the blowup of the Lie algebroid $A$ along a Lie subalgebroid $C$ supported on $N$ \cite{gualtieri.li:2014a, debord.skandalis:2021a, obster:2021a}.
More precisely, one obtains the total space of the Lie algebroid blowup by replacing the full-rank subbundle $A|_N$ by the projectivization of $((TA)|_C\setminus \ker \pi_\ast|_C)/ TC$. We denote the \emph{blowup of Lie algebroids} by the same symbol $\blup(A,C)$; 
{it is a Lie algebroid over $B$}.

The Lie algebroid structure of $\blup(A,C)$ is given by the following. Consider the space of sections of $A$ that restrict to sections of $C$, 
\begin{equation*}
    \Gamma(A,C):=\{ s\in \Gamma(A) : s|_N\in \Gamma(C)\}.
\end{equation*}
Then every $s\in \Gamma(A,C) $ canonically induces a section $\blup(s)\in \Gamma(\blup(A,C))$, and sections of this form generate $\Gamma(\blup(A,C))$, i.e.\ 
\begin{equation}\label{eq:sections_of_blup}
    \Gamma(\blup(A,C))=\mathrm{Span}_{C^{\infty}(B)} \blup(\Gamma(A,C)). 
\end{equation}
The anchor $\widetilde{\rho}$ and bracket $[\cdot,\cdot]_\blup$ are uniquely determined by
\begin{equation}\label{eq:blup_anchor_and_bracket}
    \widetilde{\rho}(\blup(s))=\widetilde{\rho(s)}\quad \text{ and }\quad [\blup(s),\blup(s')]_\blup=\blup([s,s'])
\end{equation}
for $s,s'\in \Gamma(A,C)$.

If $N$ is a transversal of the Lie algebroid $A$, another possible extension of the Lie algebroid structure to $B$ is given by the pullback Lie algebroid  $p^!A$ over $B$, which exists, since $p$ is transverse to $A$ if $N$ is. By Proposition \ref{prop:LAfol}, $p^! A$ induces the singular foliation {$p^{-1}(\cF_A)$}. 

We consider Lie subalgebroids $C\subset A$ that contain the isotropies over $N$, i.e.\ 
\begin{equation*}
    \ker(\rho|_N)\subset C.
\end{equation*}
In \S \ref{subsec:terms} we first  describe 
the singular foliation {$\cF_\blupC$} of $\blup(A,C)$ in terms of the singular foliation  {$\cF_A$} on $M$. In \S \ref{section:blup_foliation_for_transv} we assume that $C$ is supported over a transverse submanifold $N$ and express {$\cF_\blupC$} in terms of the singular foliation {$p^{-1}(\cF_A)$} of $p^! A$ on $B$.

\subsubsection{The singular foliation {$\cF_\blupC$} of $\blup(A,C)$ in terms of a singular foliation on $M$}\label{subsec:terms}

We can express {$\cF_\blupC$} of $\blup(A,C)$ in terms of {$\mathcal{F}_A$} as follows.

\begin{prop}\label{prop:blup_foliation_in_terms_of_base}
        Suppose $\ker(\rho|_N)\subset C$. Then
    \begin{equation*}
        {\cF_\blupC}=\mathrm{Span}_{C^\infty_c(B)}\{ \widetilde{Y}: Y\in {\mathcal{F}_A} \text{ such that }Y|_N\in \rho(\Gamma(C))\}.
    \end{equation*}
    Here, we denote by $\widetilde{Y}$ the unique lift of $Y$ to a vector field on $B$ which is $p$-related to $Y$.
\end{prop}

Proposition \ref{prop:blup_foliation_in_terms_of_base} follows immediately from the two following lemmas. Note that the assumption on $C$ only enters in Lemma \ref{lemma:blup_foliation_base_second_half}.

\begin{lemma}
    Let $C\subset A$ be a Lie subalgebroid. Then 
    \begin{equation*}
        {\cF_\blupC}=\mathrm{Span}_{C^\infty_c(B)}\{ \widetilde{\rho(s)} : s\in \Gamma_c(A,C) \}.
    \end{equation*}
\end{lemma}
\begin{proof}
    We have
    \begin{equation*}
        \begin{aligned}
            {\cF_\blupC}=& \widetilde{\rho}(\mathrm{Span}_{C_c^\infty(B)}\blup(\Gamma(A,C)))\\
            =&\mathrm{Span}_{C_c^\infty(B)} \widetilde{\rho}(\blup(\Gamma_c(A,C)))\\
            =&\mathrm{Span}_{C_c^\infty(B)}\{ \widetilde{\rho(s)} : s\in \Gamma_c(A,C) \}, 
        \end{aligned}
    \end{equation*}
    where we used eq.\ \eqref{eq:sections_of_blup} in the first step, properness of the blowdown map in the second (i.e.\ for $s\in \Gamma(A,C)$ the support of $s$ is compact if{f} the support of $\blup(s)$ is compact), and eq.\ \eqref{eq:blup_anchor_and_bracket} in the last. 
\end{proof}

\begin{lemma}\label{lemma:blup_foliation_base_second_half}
    Let $C\subset A$ be a Lie subalgebroid over $N$ such that $\ker(\rho|_N)\subset C$. Then
    \begin{equation*}
        s\in \Gamma(A,C) \Longleftrightarrow \rho(s)|_N\in \rho(\Gamma(C)).
    \end{equation*}
\end{lemma}
\begin{proof}
    The implication ``$\Rightarrow$'' holds trivially, the reverse is true since we assume $\ker(\rho|_N)\subset C$.
\end{proof}

\subsubsection{The singular foliation {$\cF_\blupC$} of $\blup(A,C)$ in terms of a singular foliation on $B$}\label{section:blup_foliation_for_transv}

Now, {in addition to $\ker(\rho|_N)\subset C$,} suppose that $\iota\colon N\hookrightarrow
 M$ is a transversal of $A$.
For such Lie subalgebroids, 
one has
\begin{equation}\label{eq:isoLAonB}
\blup(A,C)=\blup(p^!A,\pi_{\mathbb{P}}^! C)
\end{equation}
where the submersion $\pi_{\mathbb{P}}\colon \mathbb{P}\to N$ is the restriction of $p$.
This is a straightforward generalization of \cite[Proposition 5.16]{schuessler:2024a}, where the case $C=\iota_N^!A$ is treated.
The singular foliation {$\cF_\blupC$} 
is  necessarily tangent to the exceptional divisor $\mathbb{P}$, hence, it is distinct from {$p^{-1}(\cF_A)$}. In particular, the Lie algebroid $\blup(A,C)$ is not isomorphic to $p^!A$. It, however, comes with a canonical Lie algebroid morphism $\blup(A,C)\to p^!A \to A$.

We can express the singular foliation {$\cF_\blupC$} in terms of the singular foliation {$p^{-1}(\mathcal{F}_A)$} of $p^! A$ as follows. 

\begin{prop}\label{prop:intersection}
Suppose $\iota\colon N\hookrightarrow
 M$ is a transversal of $A$ and $\ker(\rho|_N)\subset C$.  Then 
 \begin{equation*}
     {\cF_\blupC}= {p^{-1}(\mathcal{F}_A)} \cap \mathcal{E}_C, 
 \end{equation*}
 where 
\begin{equation*}
        \mathcal{E}_C:=\{X\in \vX_c(B) : X|_{\mathbb{P}}\in \pi_{\mathbb{P}}^{-1}(\rho(\Gamma_c(C))) \}.
\end{equation*}
\end{prop}
\begin{proof}
    By eq.\ \eqref{eq:isoLAonB}, we have
    \begin{equation*}
        \Gamma_c(\blup(A,C))=\Gamma_c(\blup(p^! A, \pi^!_{\mathbb{P} }C)).
    \end{equation*}
    By eq.\ \eqref{eq:sections_of_blup} (see also \cite{gualtieri.li:2014a}), and since $\mathbb{P}\subset B$ has codimension $1$ (thus $\blup(B,  \mathbb{P})\to B$ is a diffeomorphism), we obtain that actually 
    \begin{equation*}
        \Gamma_c(\blup(A,C)) = \Gamma_c( p^! A, \pi^!_{\mathbb{P} }C).
    \end{equation*}
    Hence, {$\cF_\blupC$} is given by the intersection of the singular foliation of $p^! A$ (which is {$p^{-1}(\mathcal{F}_A)$} by Proposition \ref{prop:LAfol}) with $\mathcal{E}_C$ by Lemma \ref{lemma:blup_foliation_base_second_half}, as $\ker(\rho|_N)\subset C$ implies that the kernel of the anchor of 
$p^! A$ is contained in $\pi_{\mathbb{P}}^! C$.
\end{proof}

In some cases, we can identify or replace $\mathcal{E}_C$ by more known singular foliations on $B$, as we see in the following two examples.

\begin{example}[The restricted Lie algebroid]\label{ex:restricted_LA}
    Consider $C:=\iota^! A$. Then $\ker(\rho|_N)\subset C$ is automatically fulfilled. Since $\pi^!_{\mathbb{P}}\iota^! A = \iota_{\mathbb{P}}^! p^! A$, where $\iota_{\mathbb{P}}\colon \mathbb{P}\to B$ denotes the inclusion, we obtain that 
    \begin{equation*}
        \widetilde{\rho}(\Gamma_c(\blup(A,\iota^! A) = \rho_{p^! A}(\Gamma_c(\blup(p^! A,\iota_{\mathbb{P}}^! p^! A)) =  \rho_{p^! A}(\Gamma_c(p^! A,\iota_{\mathbb{P}}^! p^! A),
    \end{equation*}
    i.e.\ {$\cF_\blupiotaA$} consists of vector fields of the singular foliation of $p^! A$ that are tangent to $\mathbb{P}$, using again Lemma \ref{lemma:blup_foliation_base_second_half}. In other words, 
    \begin{equation*}
        {\cF_\blupiotaA}={p^{-1}(\cF_A)}\cap \Gamma(T^bB).
    \end{equation*}
    Here, $T^bB$, denotes the \emph{$b$-tangent bundle} of $B$ w.r.t.\ the hypersurface $\mathbb{P}$ (its sections are the vector fields tangent to  $\mathbb{P}$).
    Note that in general, $\Gamma_c(T^b B)\supset \mathcal{E}_{\iota^! A}$ (e.g.\ if $\iota^! A = \ker (\rho|_N)$), and only the intersection with {$p^{-1}(\cF_A)$} yields the same space.
\end{example}

\begin{example}[The isotropy Lie algebroid]\label{ex:isotropy_LA}
    Assume that $C:=\ker(\rho|_N)$ has constant rank. 
    Writing out $$\pi_{\mathbb{P}}^! C=\{(a,X)\in \pi_{\mathbb{P}}^*C\times T\mathbb{P}: \rho(a)=(\pi_{\mathbb{P}})_*X\}=
\pi_{\mathbb{P}}^*C\times \ker (\pi_\mathbb{P})_*=\rho_{p^! A}^{-1}(\ker (\pi_\mathbb{P})_*),
$$
we have by Proposition \ref{prop:intersection} that the 
induced foliation   {$\cF_\blupC$} of $\blup(A,C)$ is given by
  \begin{equation*}
  {\cF_\blupC}={p^{-1}(\cF_A)}\cap \Gamma(E).
  \end{equation*}
Here, $E$ is the \emph{edge Lie algebroid} \cite{fineedge} associated to the fibration $\pi_{\mathbb{P}}$ of the hypersurface $\mathbb{P}$;  elements in $\Gamma(E)$ are vector fields on $B$ which over  $\mathbb{P}$ are tangent to the fibers, i.e.\ lie in $\ker (\pi_\mathbb{P})_*$. 
\end{example}

   \bibliographystyle{habbrv} 

 \bigskip

\end{document}